\documentclass[11pt]{article}
\RequirePackage{amsthm,amsmath,amsfonts,amssymb}
\RequirePackage[authoryear]{natbib}
\usepackage{dsfont}
\usepackage{graphicx} 
\usepackage{natbib}
\usepackage{color}
\usepackage[OT1]{fontenc} 
\usepackage{appendix}
\usepackage{algorithm}
\usepackage{algpseudocode}
\usepackage[plainpages=false,pdfpagelabels,colorlinks=true,linkcolor=blue,urlcolor=blue,citecolor=blue]{hyperref}
\usepackage{subfigure}
\usepackage{graphicx} 		 
\usepackage{subfigure}	
\usepackage{authblk}
\usepackage{geometry}
\usepackage{amsthm,amsmath,amssymb,bbm}
\usepackage{natbib}
\usepackage{multirow}
\usepackage{subfigure}
\usepackage{makecell}
\usepackage{booktabs}
\usepackage{array}
\usepackage{tabularx}
\usepackage{tabulary}
\usepackage{caption}
\usepackage{booktabs}
\usepackage{url}
\usepackage{bm}
\usepackage{wrapfig}
\usepackage{lipsum}
\usepackage{mathrsfs}
\usepackage{dsfont}
\usepackage{titling}
\usepackage{txfonts}
\usepackage{tikz}
\usepackage{etoolbox}

\newtoggle{vaos}
\toggletrue{vaos}
\newcommand{\aosversion}[2]{\iftoggle{vaos}{#1}{#2}}

\numberwithin{equation}{section}
\theoremstyle{plain}

\newtheorem{theorem}{Theorem}[section]

\newtheorem{condition}{Condition}[section]
\newtheorem{definition}{Definition}[section]
\newtheorem{remark}{Remark}[section]
\theoremstyle{definition}

\newcommand{\abs}[1]{\left\vert#1\right\vert}
\newcommand{\norm}[1]{\left\lVert#1\right\rVert}

\newcommand{\real}{\mathds{R}}

\newcommand{\dse}{\mathds{E}}

\newcommand{\ma}{\mathcal{A}}

\newcommand{\mh}{\mathcal{H}}

\newcommand \brak [1] {\left(#1 \right)}

\newcommand{\hx}{ x}
\newcommand{\hy}{y}

\newcommand{\mb}{\mathcal{B}}

\newcommand{\ot}{\overline{T}}
\newcommand{\ut}{\underline{T}}

\title{Optimal estimation of a factorizable density
using diffusion models with ReLU neural networks\thanks{{The research was supported by ONR Grants N00014-25-1-2317 and N00014-25-1-2317, and the NSF Grant DMS-2210833.}}}
\author{Jianqing Fan\thanks{Department of Operations Research and Financial Engineering, Princeton University, jqfan@princeton.edu},\,
 Yihong Gu\thanks{Department of Biomedical Informatics, Harvard University,
 yihong\_gu@hms.harvard.edu},\, 
 and \, Ximing Li\thanks{Department of Operations Research and Financial Engineering, Princeton University, xl7408@princeton.edu}
}
\date{}
 
\begin{document}
\maketitle
\begin{abstract}
This paper investigates the score-based diffusion models for density estimation when the target density admits a factorizable low-dimensional nonparametric structure. To be specific, we show that when the log density admits a $d^*$-way interaction model with $\beta$-smooth components, the vanilla diffusion model, which uses a fully connected ReLU neural network for score matching, can attain optimal $n^{-\beta/(2\beta+d^*)}$ statistical rate of convergence in total variation distance. This is, to the best of our knowledge, the first in the literature showing that diffusion models with standard configurations can adapt to the low-dimensional factorizable structures. The main challenge is that the low-dimensional factorizable structure no longer holds for most of the diffused timesteps, and it is very challenging to show that these diffused score functions can be well approximated without a significant increase in the number of network parameters. Our key insight is to demonstrate that the diffused score functions can be decomposed into a composition of either super-smooth or low-dimensional components, leading to a new approximation error analysis of ReLU neural networks with respect to the diffused score function. The rate of convergence under the 1-Wasserstein distance is also derived with a slight modification of the method.
\end{abstract}

\section{Introduction}



The recent years have witnessed the remarkable success of diffusion models in generating images \citep{song2019generative,ho2020denoising},
audios \citep{huang2023make,kong2020diffwave}, and videos \citep{ramesh2022hierarchical,ho2022video}. See \cite{chen2024opportunities} for an overview of this subject.  Their workflow is two‑stage: a forward diffusion progressively corrupts the data with noise, which creates a map from the data distribution to a known noise distribution (usually Gaussian), and a backward process reconstructs the data distribution by solving the time‑reversal stochastic differential equations (SDEs), which provides a method that computes the inverse map from the noise distribution to the data distribution.  The backward SDE requires the diffused score function -- the gradient of the log-density along the noise path -- which is learned via score matching \citep{vincent2011connection} with neural networks. Substituting this estimated score into the reverse SDE then transforms pure Gaussian noise into samples whose distribution is close to the target distribution.

Despite their empirical success, theoretical understanding is sparse. Though estimating the generic high-dimensional distributions is fundamentally difficult and will suffer from the curse of dimensionality \citep{khas1979lower, stone1982optimal}, high-dimensional data in the real world often admits a low-dimensional structure \citep{schmidt2020nonparametric, fan2024noise, sclocchi2025phase}. This paper contributes to understanding how diffusion models with a standard neural network approximation to the score function can efficiently adapt to the unknown intrinsic low-dimensional structures. Focusing on densities that factor into lower‑dimensional components, we prove no curse-of-dimensionality convergence rates in both total variation (TV) and 1‑Wasserstein distance and establish the adaptive minimax optimality under total variation for diffusion models blind to any knowledge of function forms.  This answers an important and challenging open problem on the adaptive learning of a diffusion model to local dimensional density structure using the standard neural network architecture.

\subsection{Diffusion models and the problem under study}


The theoretical foundation for the diffusion models lies in the literature on the time-reversal of stochastic differential equations \citep{anderson1982reverse,haussmann1986time}. Let $p_0$ denote the target distribution. Consider an Ornstein-Uhlenbeck (OU) process as the forward process 
\begin{equation}\label{equation:forward}
    \mathrm{d} X_t=-\beta_t X_t \mathrm{~d} t+\sqrt{2 \beta_t} \mathrm{~d} B_t, \qquad X_0\sim p_0
\end{equation}
where $\left(B_t\right)_{t\geq 0}$ is a standard Brownian motion and $\beta_t:[0,\infty) \to [0,\infty)$ is user-chosen weighting function. Denote the marginal density of $X_t$ by $p_t(\cdot)$. For a fixed time endpoint $\overline{T}$ specified later, it is well-known that the backward process $\left(Y_t\right)_{[0, \ot]}$ with $Y_t=X_{\ot-\mathrm{t}}$ is the solution to the following SDE: 
\begin{align}\label{equation:backward}
    & d Y_t=\beta_{\ot-t}\left(Y_t+2 \nabla \log p_{\ot-t}\left(Y_t\right)\right) d t+
    \sqrt{2 \beta_{\ot-t}} d Z_t, \qquad Y_0 \sim p_{\ot},
\end{align}
with $\left(Z_t\right)_{t\geq 0}$  is another standard Brownian motion and $\nabla \log p_t$ is referred to as the score function for $X_t$.  Because the OU process rapidly drives $p_{\ot}$ toward the standard Gaussian as $\ot \to \infty$ \citep{bakry2014analysis},  initializing $Y_0$ from Gaussian noise and running model \eqref{equation:backward} yields a sample exactly from $p_0$ -- provided the score is known. In practice, one discretizes the SDE and learns an estimated $\nabla_x\log p_t(x)$ from data, thereby obtaining an implementable diffusion model.  


The \emph{diffusion model}, motivated by the above SDE, uses the following procedure to estimate the target distribution implicitly. We first specify the forward process \eqref{equation:forward} that determines how the target $p_0$ is transformed to Gaussian noise. Given $n$ observations sampled from $p_0$, we can apply the score-matching to estimate the score functions $s(x,t) = \nabla_x \log p_t(x)$, typically using neural networks. Finally, plugging in the estimated score function $\hat{s}$ into the reverse process \eqref{equation:backward} allows it to generate samples $\widehat{X}$ whose distribution is close to $X\sim p_0$ from the noise.  For an overview, see \cite{tang2025score}.


There is a considerable literature on offering statistical analyses of the aforementioned diffusion model, assuming $p_0$ belongs to a certain function class; see, for example, \cite{chen2023score, oko2023diffusion}. Focusing on diffusion model using fully connected neural networks for score function estimation, \cite{oko2023diffusion} first shows it can attain minimax optimal error rate both in TV and 1-Wasserstein distance when $p_0$ belongs to H{\"o}lder class. The key idea is that neural networks can approximate the diffused score function in a nearly optimal way under this scenario. Given the key part in the paradigm is to estimate the (diffused) score function from data, there is also a considerable literature focusing on understanding the optimal $L_2$ estimation error of the score function \citep{wibisono2024optimal, zhang2024minimax, dou2024optimal} with more traditional nonparametric techniques like the kernel density estimator. However, all these convergence rates suffer from the curse of dimensionality that results in an extremely slow rate when the dimension $d$ is large.


Real‑world data -- ranging from physical laws \citep{dahmen2022compositional} to images and language \citep{partee1984compositionality} -- often exhibit certain low-dimensional compositional structures. Neural networks are adept at exploiting this and can algorithmically learn low‑dimensional compositional structures without assuming explicit function form \citep{bauer2019deep, schmidt2020nonparametric,kohler2021rate, fan2024factor}.
Diffusion models appear to share the same strength: empirical studies indicate they capture the low-dimensional patterns \citep{sclocchi2024probing,sclocchi2025phase}, and recent theoretical studies show that the backward sampling dynamic \eqref{equation:backward} is provably adaptive to low-dimensionality in a structure-agnostic way \citep{potaptchik2024linear, huang2024denoising, li2024adapting} in terms of discretization error for solving \eqref{equation:backward} with a known score function. From a statistical viewpoint, however, key questions remain unanswered: It still lacks a clear understanding of whether the entire diffusion model pipeline, including estimating the score function, is sample efficient under the unknown low-dimensional structure. In particular, the sample efficiency of the score‑function estimation stage has not been fully characterized. Although there is some work showing diffusion models \citep{tang2024adaptivity, azangulov2024convergence} can adapt to the unknown manifold, these are essentially similar to the H{\"o}lder case studied by \cite{oko2023diffusion}, where the density can be expanded as a series of fixed local polynomial bases and minimax optimal rate can be attained via a variant of kernel density estimation \citep{tang2023minimax}. 


In this paper, we study whether diffusion models can efficiently estimate distributions from the exponential‑interaction family, a subclass of the general hierarchical composition structure. Although this family exhibits a seemingly simple low‑dimensional structure, it falls far beyond the settings explored by current diffusion model theory, and by taking interaction order sufficiently large, it can include most of the frequently used low-dimensional models.  However, in the statistical machine learning community, there is still no practical and sample-efficient estimator for such a density, along with a corresponding algorithm to generate the data from the estimated density. 
To be specific, let $d$ be the ambient dimension and suppose we observe i.i.d. samples $\{X_{0, i}\}_{i=1}^n$ from an unknown density $p_0$ that admits   the form
\begin{equation}\label{equation:interactionintro}
    p_0(x)=\exp \left(\sum_{J \in \mathcal{S}} f_J(x_J) \right) \qquad \text{with} \qquad \mathcal{S}\subseteq \{J\subseteq \{1,\ldots, d\}: |J|= d^*\}
\end{equation} where $x_J=(x_j)_{j\in J}$ and each $f_J$ is a $d^*$-variate $\beta$-H{\"o}lder function. Our central question is:
 \begin{center}
    \it
    Can diffusion models with fully connected neural networks as score function estimators circumvent the curse of dimensionality under \eqref{equation:interactionintro} in a structure-agnostic manner?
 \end{center}
This question is very technically challenging, as the score functions at most time points do not admit a low-dimensional structure.  Dedicated efforts are made to answer this critical but challenging question in generative AI. 

The density structure \eqref{equation:interactionintro} depicts the functional form of the distributions that appear in many graphical model frameworks, such as Markov random fields and Bayesian networks. These models succinctly encode local dependencies in high‑dimensional data \citep{besag1974spatial,cross1983markov,saul1994boltzmann} and are therefore widely adopted in modern statistical practice \citep{wainwright2008graphical}. Despite this popularity in statistical modeling, it still lacks an estimator that is both theoretically minimax‑optimal and easily implementable. Existing guarantees are largely confined to special cases such as tree‑structured graphs \citep{liu2011forest, gyorfi2023tree}. In more general settings, \cite{vandermeulen2024breaking, vandermeulen2024dimension} proposes explicit density estimators $\hat{p}$, but sampling from them is computationally expensive, and the convergence rate is not optimal. \cite{kwon2025nonparametric} brings diffusion models into the picture, but they adopted their specially designed neural network that essentially performs integration, which is impractical in implementation. In short, there is still no simultaneously practical and provably sample‑efficient procedure for learning and generating samples from densities of the form \eqref{equation:interactionintro}.

Unlike previous manifold settings that are either linear \citep{oko2023diffusion, chen2023score} or nonparametric \citep{tang2024adaptivity, azangulov2024convergence} whose low-dimensional structure of similar form is maintained throughout all the timesteps $t$, the main challenge here is that this low-dimensional factorization structure no longer holds for diffused $p_t$ when $t$ is not very small.
From a technical perspective, even if the interaction set $\mathcal{S}$ were known, one cannot easily construct a fixed set of basis functions that can approximate densities of the form \eqref{equation:interactionintro} with a dimension‑independent error, indicating that one cannot get rid of the curse of dimensionality via previous proof strategies in \cite{oko2023diffusion,tang2024adaptivity,azangulov2024convergence}. 
Addressing this challenge calls for new insights and tools, which are also essential for the understanding of diffusion models under the general hierarchical composition structures.

\subsection{Contributions}

This paper answers the question in the affirmative. To be specific, we show that under the general regularity condition that the density $p_0$ is supported on a unit cube and is bounded from below and above, the diffusion model can obtain the error rate $n^{-\frac{\beta}{2\beta + d^*}}$ in TV distance if $p_0$ admits the low-dimensional factorized structure \eqref{equation:interactionintro}. The diffusion model we analyzed employs a standard fully connected ReLU neural network that estimates the score function for all $t$, without any prior knowledge of the underlying factorization. This matches the minimax optimal convergence rate in TV distance when $p_0$ lies within this exponential interaction family and has implications both in diffusion model theory and structured density estimation and sampling.
\begin{itemize}
\item We prove that a vanilla diffusion model -- equipped with a fully connected ReLU network for score estimation -- adaptively exploits the low‑dimensional factorization in \eqref{equation:interactionintro}. This composition of a super‑smooth link and an interaction term departs sharply from previously analyzed settings and unveils, for the first time, how diffusion models learn such a hierarchical structure.
\item Simulating the reverse SDE to TV distance $\epsilon$ only requires $O(\mathrm{poly}(\log(d),d^*)\cdot  \epsilon^{-1})$ timesteps \citep{li2024adapting, huang2024denoising}. Coupled with our statistically minimax-optimal rate of convergence, this makes the vanilla diffusion model, to the best of our knowledge, the first method that is simultaneously provably optimal and readily implementable for sampling from the general nonparametric family in \eqref{equation:interactionintro}.
\end{itemize}

Furthermore, following the idea of piecewise score estimator in \cite{oko2023diffusion}, we show that one can also attain estimation error $n^{-\frac{\beta+d^*/d}{2\beta+d^*}}=o(n^{-\frac{\beta}{2\beta + d^*}})$ under 1-Wasserstein distance. 

The main technical novelty is to derive the dimension-free approximation error of the score function $s_t(x) = \nabla_x \log p_t(x)$ across $t$. 
Unlike the previous analysis, which utilizes the fact that $p_0$ itself can be decomposed into a series of bases with an optimal approximation error, we instead decompose the diffused score $s_t(x)$ into a composition of either super-smooth \citep{fan1991optimal} or low-dimensional functions throughout $t$. Formally, let $\mathsf{c}>0$ be any fixed large constant. We show that, for any given depth $L$ and width $W$, there is a neural network $\tilde{f}$ satisfying
\begin{align*} 
   \forall t \in [(WL)^{-\mathsf{c}}, \mathsf{c}\cdot \log(WL)] \qquad \left[\int_{x\in \real^d}\norm{\tilde{s}(x, t) - s_t(x)}_2^2 p_t(x)dx\right]^{1/2} \lesssim \frac{(WL)^{-2\beta/d^*}}{\sqrt{t}\wedge 1}.
\end{align*}

Besides circumventing the curse of dimensionality, our approximation error and stochastic error analyses sharpen prior results \citep{oko2023diffusion,kwon2025nonparametric} even when $d=d^*$. Particularly, we get rid of the boundary super‑smoothness assumption, require only density bounded from below and above, and obtain explicit error rates for arbitrary network widths and depths without imposing sparse weights. Meanwhile, we establish a stochastic error bound when the weights of neural networks are unbounded or grow exponentially with $(L, W)$, while the previous stochastic analysis inevitably requires polynomially bounded weights to control the $L_\infty$ cover number.

\subsection{Related works and comparisons}

Underlying the success of modern deep learning models is the approximation capability of neural networks. Classical universal-approximation results show that a one-hidden-layer neural network can approximate any continuous function on a compact domain to arbitrary precision \citep{hornik1991approximation,barron1993universal}. Yet universality alone is not distinctive -- splines \citep{friedman1991multivariate}, wavelets \citep{donoho1995adapting}, and other non-parametric techniques enjoy the same property. Recently, \cite{telgarsky2016benefits} illustrated the benefits of deep neural networks over shallow ones. Based on the the idea of approximating the local polynomial expansion of the target function and the fact that polynomial functions can be efficiently approximated by the constructed neural networks in \cite{telgarsky2016benefits} and \cite{yarotsky2017error}, there is a considerable literature on quantifying how the $L_2$ or $L_\infty$ error of a fully connected ReLU network scales with its depth $L$ and width $W$ when the target functions lies within some specific function classes \citep{yarotsky2017error, shen2022optimal, kohler2021rate, lu2021deep}. However, these rates suffer from the curse of dimensionality in that they will be slow even for moderate $d$. \cite{schmidt2020nonparametric, kohler2021rate} show that neural networks can circumvent the curse of dimensionality when the target function is a composition of simple functions. And there is also some literature characterizing the provable advantages of neural networks over traditional nonparametric techniques via establishing the lower bounds for the latter, like wavelets \citep{schmidt2020nonparametric} and general linear estimators \citep{suzuki2021deep, imaizumi2022advantage}. In the context of the diffusion model, where the function of interest is the diffused score function $f^\star(x, t)=\nabla_x \log p_t(x)$,  \cite{oko2023diffusion} developed the approximation ability of neural networks to the diffused score function using the fact that the target density $p_0$ can be approximated by a sequence of fixed basis functions at the optimal rate -- a strategy later extended to unknown manifolds \citep{tang2024adaptivity,azangulov2024convergence}. 
However, one cannot simply apply similar arguments to the factorized structure we studied.
We propose a novel decomposition of the diffused score function, showing that a low-dimensional compositional structure remains in $f^\star$ even when $t$ is not very small, the regime where the factorized structure breaks down; see the discussion in Section \ref{section:approximationtheorem}.

There is a considerable literature on the theoretical understanding of diffusion models. When an accurate score function is provided, reverse-time SDEs generate samples that are close to the target distribution in polynomial time \citep{chen2022sampling,de2022convergence,lee2022convergence,li2024d}, and the iteration complexity remains adaptive to the underlying low-dimensional structure \citep{li2024adapting,huang2024denoising}. In the statistical setting where the score should be learned from data, diffusion models trained by score matching achieve minimax-optimal rates for H{\"o}lder-smooth densities \citep{oko2023diffusion}, for data supported on unknown manifolds \citep{tang2023minimax,azangulov2024convergence}, and for certain parametric families \citep{mei2025deep}. Because classical approaches like kernel methods \citep{tang2023minimax} attain the same optimal rates in these regimes, the distinct advantage of diffusion models remains ambiguous. \cite{cole2024score} showed that a vanilla diffusion model can break the curse of dimensionality in a structure-agnostic way, but their bounds require the log-density to be nearly Gaussian and fail to attain optimal rates. \cite{kwon2025nonparametric} considered the exponential-interaction family studied here, but leveraging an impractical, hand-crafted neural architecture that integrates the density -- an architecture substantially different from the neural networks commonly used in practice. Therefore, it is still unclear whether the vanilla diffusion model adopted in practice can learn low-dimensional composition structures that neural networks are adept at. This paper makes progress in this direction. 

Beyond diffusion models, deep generative alternatives -- generative adversarial networks (GANs), variational autoencoders (VAEs), and normalizing flows -- have been analyzed for their statistical efficiency. Recent studies (e.g., \citealp{liang2021well, singh2018nonparametric, uppal2019nonparametric, huang2022error, chen2020distribution, arthur2024wgan}) provide finite-sample guarantees, but almost all assume the target density belongs to a function class where classical methods can also obtain the optimal rate. Another line of work \citep{bos2024supervised, vandermeulen2024breaking, vandermeulen2024dimension} proposes new algorithms that can obtain rates that beat the classical curse-of-dimensionality bound $n^{-\beta/(2\beta+d)}$ under \eqref{equation:interactionintro}. However, these rates are still sub-optimal, and the resulting estimators do not furnish a practical procedure for drawing samples -- two limitations that the vanilla diffusion model we analyzed here overcomes.

\subsection{Organization}
The rest of the paper is as follows. 
In Section \ref{section:setups}, we introduce the setup, in particular, the diffusion model and score-matching. 
In Section \ref{section:mainresult}, we first provide the approximation error for the score function using a deep fully-connected ReLU network with any architecture hyperparameters. 
Based on this generic approximation result, we establish the convergence rate for the distribution generated by the diffusion model (with optimally tuned network parameters). 
Section \ref{sec:sketch} briefly sketches the proof for the main results and presents our novel technical tools for 
deep score approximation. 
Section \ref{section: w1} provides the modification of the diffusion model and the corresponding estimation rate in $W_1$ distance. 
All detailed proofs are contained in the supplemental material and are available upon request.
 
\subsection{Notations}

The following notations will be used throughout this paper. We use $c_1, c_2, \ldots$ to denote the global constants that appear in the statement of any theorem, proposition, corollary, and lemma. We use $C_1, C_2, \ldots$ to denote the local intermediate constants in the proof. Hence, all the $c_1, c_2, \cdots$ have unique referred numbers, while all the $C_1, C_2, \ldots$ will have different referred numbers in respective proofs. 
The notation $C(A_1,A_2,\cdots,A_n)$ means that the constant $C$ only relies on $A_1,A_2,\cdots,A_n.$
We use $a(n) \lesssim b(n)$, $b(n) \gtrsim a(n)$, or $a(n) = O(b(n))$ if there exists some constant $C>0$ such that $a(n) \le Cb(n)$ for any $n \ge 3$, We use $\widetilde{O}(b(n))$ to hidden poly-log factors. Denote $a(n) \asymp b(n)$ if $a(n)\lesssim b(n)$ and $a(n) \gtrsim b(n)$. 
For a vector $x$, we use $\norm{x}_p$ to denote the $\ell_p$-norm of $x$ for $0<p\leq +\infty$, and 
the $\ell_2$-norm of $x$ will be abbreviated as $\norm{x}$. 
For any $m\in \mathds{N}^+$, we use $[m]$ to denote the set $\{1,2,\cdots,m\}$. 

\section{Setup and methodology}\label{section:setups}

Consider the following implicit distribution estimation problem. Let $p_0$ be the unknown $d$-variate density function, and suppose we observe $n$ i.i.d. $d$-dimension samples $ \{X_{0, i}\}_{i=1}^n\ $ drawn from $p_0$. The goal of implicit distribution estimation is to find a transform $\hat{h}$ of the noise $U$ such that the distribution $\widehat{X}=\hat{h}(U)$ is close to $p_0$. Note that when $U$ is the uniform distribution, $\hat h$ estimate the inverse cumulative distribution function of $p_0$.  Let $\hat{p} = p_{\widehat{X}}(x)$ be the density function of the random variable $\widehat{X}$ conditioned on the data (with randomness $U$).  One can use the TV distance between $\hat{p}$ and $p_0$, defined as 
\begin{align*}
    \mathsf{TV}(\hat{p}, p_0) = \int |\hat{p}(x) - p_0(x)| dx,
\end{align*} 
to evaluate the performance of the implicit distribution estimator ${p}_{\widehat X}(\cdot)$ or the quality of the simulated sample $\widehat{X}$ from the target distribution $p_0$.

We first introduce some notations used in the construction of the diffusion model estimator. Let $\beta_t: \mathds{R}^+\to \mathds{R}^+$ be a pre-determined weighting function. It is easy to see that  the solution to the stochastic differential equation \eqref{equation:forward} is given by
$$
X_t = m_t X_0 + m_t \int_0^t \sqrt{2 \beta_s} / m_s dB_s,
$$
where $m_t=\exp(-\int_0^t \beta_s ds)$.
Hence,  the conditional distribution $X_t|X_0$ in \eqref{equation:forward} is a Gaussian $\mathcal{N}(m_tX_0,\sigma^2_t)$ where
$\sigma_t^2=1-m_t^2$.  
Thus the conditional distribution $X_t|X_0$ has the density
\begin{align}\label{equation:condition}
    p_{t|0}(x|X_0=x_0) = \frac{1}{\sigma_t^d (2\pi)^{d/2}}\exp\left(-\frac{\norm{x-m_t x_0 }^2}{2\sigma_t^2}\right).
\end{align}
We denote the marginal density $p_t(x)=p_{X_t}(x)=\int p_{t|0}(x|y) p_0(y) dy$, and the score function $s_t(x) = \nabla_x \log p_t(x)$. We suppress the dependencies of $p_t$ and $s_t$ on $\beta_t$ for simplicity. The diffusion model first estimates the score function $s_t(x)$ by de-noising score matching, and then plugs the estimated score function $\hat{s}_t(x)$ into the reverse SDE in \eqref{equation:backward} to construct the estimator $\widehat{X}$, which is a realization from $\hat p$ with an associated algorithm by discretizing \eqref{equation:backward}.

In this article, we analyze the performance of diffusion models with a standard configuration that performs the score matching using fully-connected deep ReLU neural networks. We show that the proposed estimator adapts to the low-dimensional structure of $p_0$ in a structure-agnostic manner. We first introduce the neural network classes, score-matching estimator, and the constructed diffusion model in Section \ref{section:algorithm}, and present the considered low-dimension structure in Section \ref{subsec:interaction-model}.

\subsection{The proposed diffusion model estimator}\label{section:algorithm}

\noindent \textbf{Fully connected deep ReLU neural networks.} We adopt the fully connected deep neural network with ReLU activation $\sigma(\cdot) = \max\{\cdot, 0\}$ to estimate the score function given its great empirical success. We refer to it as \emph{deep ReLU network} for short. Let $L$ and $W$ be any positive integers.  A deep ReLU network with depth $L$ and width $W$ admits the form of
\begin{align}
\label{equation:nn-def}
f(x)=\mathcal{T}_{L+1} \circ \bar{\sigma}  \circ  \mathcal{T}_L  \circ  \bar{\sigma}  \circ  \mathcal{T}_{L-1}
\circ  \bar{\sigma}  \circ  \cdots   \circ \mathcal{T}_2  \circ \bar{\sigma} \circ  \mathcal{T}_1(\hx),
\end{align}
Here $\mathcal{T}_i(\hx)=A_i \hx+b_i$ is a linear transformation with $A_i \in \real^{d_i \times d_{i-1}}, b_i \in \real^{d_i}$, where the natural number sequence $\left(d_0, d_1, \cdots, d_L, d_{L+1}\right)$ is specified as $(d, W, \cdots, W, 1)$, and $\bar{\sigma}: \mathds{R}^{d_i} \rightarrow \real^{d_i}$ applies the ReLU activation function $\sigma(x)=\max \{0, x\}$ to each entry of a ${d_i}$-dimensional vector. We define the family of \emph{deep ReLU networks taking $d$-dimensional vector as input with depth $L$ and width $W$} as
\begin{align}
\label{equation:nn-without-trunc}
\mathcal{F}_{\mathsf{NN}}(d, L, W) = \{g(x): g(x) ~\text{in}~ \eqref{equation:nn-def}\}.
\end{align}


\medskip 
\noindent \textbf{Score matching.} 
It is required to first estimate the \textit{score function} $s_t(x)=\nabla \log p_t(x)$ 
using observed data $\{X_{0,i}\}_{i=1}^n$ and then plug-in the estimated $\hat{s}_t(x)$ into the backward process \eqref{equation:backward} for sampling. A standard method to estimate the score function is the \emph{de-noising score-matching} \citep{hyvarinen2005estimation,vincent2011connection}.

Recall the definition of $p_{t|0}$ in \eqref{equation:condition}. \cite{vincent2011connection} establishes the following identity, which express the mean-square error for the score estimation as that for the conditional one.  That is, for any $s: \mathds{R}^{d+1} \to \mathds{R}^d$,
\begin{align}\label{equation:scoreidentity}
\begin{split}
   &\dse_{X_t}\left[\norm{s(X_t,t)-\nabla \log p_t(X_t)}_2^2\right] \\
&   \qquad \qquad =\dse_{ X_t\sim p_{t|0}(\cdot|X_{0}),X_0}\left[\norm{s(X_t,t)-\nabla \log p_{t|0}(X_t\mid X_0)}_2^2\right]+C_{t,p_0},
\end{split}
\end{align}
where $C_{t,p_0}$ is a constant only relies on $t$ and $p_0$ but independent of $s$.  The beauty of this identity is that there is an explicit form for the conditional density function, and the constant is independent of $s$, and hence the neural network weights in our approximation.

Let $(\underline{T}, \overline{T}) \in \mathds{R}^2$ be time hyper-parameters that $0 < \underline{T} < \overline{T}$. 
Motivated by \eqref{equation:scoreidentity}, for any $s: \mathds{R}^{d+1} \to \mathds{R}^d$, the finite sample score matching loss over time interval $[\underline{T}, \overline{T}]$ is defined as
\begin{align}\label{equation:loss}
    \widehat{\mathsf{L}}_n(s) = \frac{1}{n} \sum_{i=1}^n \mathds{E}_{t\sim \mathsf{Unif}(\underline{T}, \overline{T}), X_t\sim p_{t|0}(\cdot|X_{0,i})} \left[\lambda(t)\left\|s(X_t, t)-\nabla \log p_{t|0}(X_t | X_{0, i})\right\|_2^2\right],
\end{align} 
where $\{X_{i, 0}\}_{i=1}^n$ are given data, $\mathsf{Unif}(\underline{T},\overline{T})$ is the uniform distribution  on the time interval $[\underline{T},\overline{T}]$. Note that $p_{t|0}$ has a parametric form in \eqref{equation:condition} such that one can easily sample $X_t$ and $\nabla_x \log p_{t|0}(x|z) = - (x-m_t z)/\sigma_t^2$. In this paper, we choose $\lambda(t)\equiv 1$.

\begin{remark}
Note that the loss function in \eqref{equation:loss} involves continuous time sampling and expectation over the normal distribution $p_{t|0}$. In practice, the expectation over $T\sim \mathsf{Unif}(\underline{T},\overline{T})$ and $X_t \sim p_t|0(\cdot|X_{0,i})$ can be approximated by Monte Carlo sampling and further be optimized using gradient descent algorithm \citep{kloeden2011multilevel,song2019generative,de2022convergence,benton2023linear}. In this paper, we focus on the semi-empirical loss \eqref{equation:loss} and solve the continuous-time SDE in the next step for simplicity of presentation. We emphasize that our main result can be easily extended to the fully empirical level loss with a proper time discretization by adopting a similar idea in Section 5.3 in \cite{oko2023diffusion}. 
\end{remark}

The neural network score-matching estimator is 
\begin{align}\label{equation:estimator}
\hat{s}:=\arg\min_{s\in \mathcal{F}_{\mathsf{NN}}(d+1,L,W,B)}
\widehat{\mathsf{L}}_n(s),
\end{align} where the truncated neural network class is defined as
\begin{align}\label{eq:nn-truncationdef}
    \mathcal{F}_{\mathsf{NN}}(d+1, L, W, B) =
    \left\{\tau\Big(f(x, t); B \sigma_t^{-1}\sqrt{\log (WL)}\Big), f(x, t)\in \mathcal{F}_{\mathsf{NN}}(d+1, L, W)\right\}
\end{align} where $\tau(z; \rho) = \mathrm{sign}(z) \min\{|z|, \rho\}$. Here $B$ is the truncation hyper-parameter to be determined.
 
\medskip
\noindent \textbf{Diffusion model with de-noising score-matching.} 
We propose the distribution estimator $\widehat{X}$ for $X$ by solving the backward process 
with the estimated score function $\hat{s}$ in \eqref{equation:estimator}. More precisely, we regard the density of $\widehat{X}$, conditioned on the data, as an implicit density estimator of $p_0$. Additionally, equation \eqref{equation:backwardestimated} also provides a method to sample the data $\widehat{X}$.  Since $p_{\overline{T}}$ in \eqref{equation:backward} converges to the standard Gaussian in an exponential rate as $\overline{T}\rightarrow+\infty$ \citep{bakry2014analysis}, we pick a large enough $\overline{T} \asymp \log(n)$ and replace $p_{\overline{T}}$ in \eqref{equation:backward} by standard normal distribution. The process is summarized in Algorithm \ref{algorithm: diffusion}. Note that $\widehat{X}$ is a transform of the noise $\widehat{Y}_0$ and Brownian motion $\{B_t\}_{t\ge 0}$, which serves as an implicit estimator for the target distribution $X$.


\begin{algorithm}[!t]
    \renewcommand{\algorithmicrequire}{\textbf{Input:}}
    \renewcommand{\algorithmicensure}{\textbf{Output:}}
    \caption{Implicit distribution estimator based on diffusion models}
    \label{algorithm: diffusion}
    \begin{algorithmic}[1]
        \Require i.i.d. samples $\{X_{0,i}\}_{i=1}^n$ drawn from the distribution $p_0$.
        \Ensure  an output $\widehat{X}$ drawn for an implicitly estimated density that is close to $p_0$ for large $n$ 
        \State Set the network parameters (depth $L$, width $W$, and truncation threshold $B$) 
        and time parameters $(\underline{T},\overline{T})$ as specified in Condition \ref{condition:parameter}. 
        \State Compute the estimated score function
        $\hat{s}$ by solving $\hat{s}:=\arg\min_{s\in \mathcal{F}_{\mathsf{NN}}(d+1,L,W,B)} \widehat{\mathsf{L}}_n(s)$, 
        where the loss $ \widehat{\mathsf{L}}_n$ is defined in \eqref{equation:loss}.
        \State Draw the sample from the backward process with the estimated score $\hat{s}$:
        \begin{align}\label{equation:backwardestimated}
        & d \widehat{Y}_t=\beta_{\overline{T}-t}\left(\widehat{Y}_t+2 \hat{s}( \widehat{Y}_t,\overline{T}-t)\right) d t+
        \sqrt{2 \beta_{\overline{T}-t}} d B_t, \, \widehat{Y}_0 \sim N(0,I_d).
        \end{align}
    \State \Return $\widehat{X}$ as $\widehat{Y}_{\overline{T}-\underline{T}}$.
 \end{algorithmic}
\end{algorithm}

\subsection{Exponential-interaction model}
\label{subsec:interaction-model}

We first introduce the $(\beta, C)$-smooth function class.  
\begin{definition}[$(\beta,C)$-smooth]\label{definition:holder}
Consider $\beta=q+s$ for some natural number $q$ and $0<s \leq 1$, and $C>0$.
A $d$-variate function $f:\real^d\rightarrow\real$ is $(\beta, C)$-smooth if for every sequence $\left\{\alpha_j\right\}_{j=1}^d$ of
natural numbers such that $\sum_{j=1}^d \alpha_j=q$, the partial derivative $(\partial f) /\left(\partial x_1^{\alpha_1} \cdots \partial^q x_d^{\alpha_d}\right)$ exists and satisfies for any $x, y \in \mathds{R}^d$ that
$$
\left|\frac{\partial^q f}{\partial x_1^{\alpha_1} \cdots \partial x_d^{\alpha_d}}(\hx)-\frac{\partial^q f}{\partial y_1^{\alpha_1} \cdots \partial y_d^{\alpha_d}}(\hy)\right| \leq C\|\hx-\hy\|_2^s .
$$
\end{definition}

It is known that estimating the density for the $d$-variate $(\beta,C)$-function class 
from $n$ i.i.d. observations has the minimax estimation error as $n^{-2\beta/(2\beta+d)}$ \citep{khas1979lower,stone1982optimal}. The rates are extremely slow when the input dimension $d$ is very large, with the effect degree of smooth $\beta/d$ \citep{fan2024noise}, which is common in neural network applications \citep{poggio2017and}. In this paper, we consider the following low-dimensional structure.

\begin{definition}[Exponential-interaction model]\label{definition:interaction}
Given a nature number $d^*\leq d$ and real numbers $\beta,C>0$, the exponential-interaction function class 
$\mathcal{H}(d,d^*,\beta,C)$ is defined as follows
\begin{align*}
    \mathcal{H}(d, d^*, \beta, C) = \left\{p(x)=\prod_{J\subseteq [d], |J|\le d^*} f_J(x_{J}): ~f_J \text{ is }(\beta,C)\text{-smooth}\right\}.
\end{align*}
\end{definition}

\begin{remark}\label{remark:example}
A similar model class (additive rather than productive) has been studied in \cite{bhattacharya2024deep} for understanding the performance of high-dimensional neural network regression.  We provide several examples to show the versatility of our model class. 

\medskip
\noindent \textbf{(a)}. 
Our model class $\mh(d,d^*,\beta,C)$ captures any $(\beta,C)$-smooth density $p:\real^{d^*}\rightarrow\real$.
For any $J\subset [d]$ with size $d^*,$ we set 
$f_J= p$ if $J=[d^*]$ and $1$ else wise. 
In this way, Definition \ref{definition:interaction} recovers the density $p$. 

\medskip
\noindent\textbf{(b)}.
By the Hammersley-Clifford theorem \citep{wainwright2019high},  our distribution class  
  captures any distribution satisfying the local Markov property with respect to a graph with $d$ vertices whose maximum clique size is $d^*$. For instance,  the Markov random field in the left panel of Figure \ref{fig:example} belongs to our density class with $d=4$ and $d^*=2$, which has the density $ p(x) \propto f_{12}(x_1,x_2)f_{23}(x_2,x_3)f_{34}(x_3,x_4)f_{14}(x_1,x_4) $.

\medskip

\noindent\textbf{(c)}. Our distribution class also covers Bayesian networks. 
For instance,  the Bayesian network in the right panel of Figure \ref{fig:example}
  belongs to our density class with $d=4$ and $d^*=2$, which has the density 
    $ p(x) \propto p(x_1)p(x_2\mid x_1)p(x_3\mid x_2)p(x_4\mid x_3) $. 
\begin{figure}[t!]
    \centering
    \subfigure[Markov random field]{\includegraphics[width=0.25\linewidth]{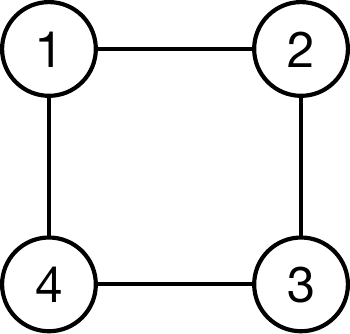}}
     \qquad
      \qquad
    \subfigure[Bayesian network]{\includegraphics[width=0.25\linewidth]{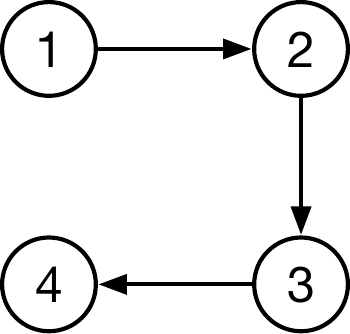}}
    \caption{The illustrating examples with 4 variables ($d=4$) and $d^*=2$. For the Markov random field, the edge between node $x$ and $y$ 
    indicates that variables $x$ and $y$ may not be independent given all the other variables. 
    For the Bayesian network, the arrow from node $x$ to $y$ 
    indicates that $x$ is the direct cause of $y$.}
        \label{fig:example}
\end{figure}
 
\end{remark}
%


\section{Main result}\label{section:mainresult}

We first impose some regularity conditions widely adopted in the density estimation literature. 
\begin{condition}[Lower and upper bounded]\label{condition:lowerbound}
The density $p_0(x)$ is supported on the unite cube $[-1,1]^d$. There exists a constant $c_1>1$ such that $1/c_1\leq p_0(x)\leq c_1$ for any $x\in [-1,1]^d$.
\end{condition}

Note that the diffusion model estimator is based on the weighting function $\beta_t$ since the target score function $s_t=\nabla \log p_t(x)$ and the reverse backward process \eqref{equation:backward} depend explicitly on $\beta_t$. Condition \ref{condition:betat} below shows that our result applies to various choices of $\beta_t$ satisfying some weak regularity conditions. This includes some widely adopted choices of $\beta_t$, including, for example, the constant function in the Langevin dynamics \citep{song2019generative} and the linear function in the DDPM \citep{ho2020denoising}, in previous literature.

\begin{condition}\label{condition:betat}
There exists a constant $c_2>1$ such that for any $\ell\in \mathds{N}$ and $t>0$, the derivatives $\abs{\frac{d^{\ell}}{dt}\beta_t}\leq c_2$ and $1/c_2\leq \beta_t\leq c_2$.
\end{condition}

\subsection{Score approximation error}\label{section:approximationtheorem}

Recall $s_t(x) = \nabla \log p(x)$ is the score function of the diffused density $p_t$ in Section \ref{section:setups}. The following theorem establishes the approximation error of fully connected neural networks with respect to the function $s^*(t, x)=s_t(x): \mathds{R}^{d+1}\to \mathds{R}^d$ when $p_0 \in \mathcal{H}(d, d^*, \beta, C)$.  The technical proof is highly complex and technical, with over 100 pages.  We therefore relegate all the proofs to the appendix.

\begin{theorem}[Approximation error]\label{theorem:approximation}
Assume Conditions  \ref{condition:lowerbound} and \ref{condition:betat} hold.  
There exists constants $c_3,c_4,c_5, c_6,c_7$ depending  only on $(c_1,c_2, \beta,d,d^*,C)$ such that the following holds. 
For any $W,L$ satisfying $\min\{W, L\}\geq \brak{1+\log (WL)}^{c_5}$, letting $\ut=(WL)^{-c_3}$ and $\ot=c_4\log (WL)$,  and any $p_0 \in \mh(d,d^*,\beta,C)$,
there exists a neural network $\tilde{s} \in \mathcal{F}_{\mathsf{NN}}(d+1, L, W, c_6)$ such that 
\begin{align}
\label{equation:approx-error}
      \left[\int_{x\in \real^d}\norm{\tilde{s}(x, t) - s_t(x)}^2 p_t(x)dx\right]^{1/2} \le c_5 \frac{(WL)^{-2\beta/d^*}}{\sigma_t}
      \big(\log (WL)\big)^{c_7} \qquad \forall t\in [\ut,\ot].
\end{align} The dependency of $c_3,\ldots,c_6$ on $(c_1,c_2, \beta,d,d^*,C)$ can be found in 
\aosversion{Condition D.1}{Condition \ref{condition:parameterproof}}.
\end{theorem}

Theorem \ref{theorem:approximation} is the first neural network score function approximation result that can further yield the optimal estimation error rate when $\log p_0(x)$ admits a low-dimensional interaction model structure.  It indicates that the approximate error depends on the effective dimension-adjusted degree of smoothness $\beta/d^*$\citep{fan2024noise}.  There is a considerable literature on establishing deep ReLU networks' approximation ability on the diffused score function $s_t(x)$ when $p_0$ belongs to specific function classes. For example, \cite{oko2023diffusion} shows if $p_0$ can be approximated by a sequence of fixed basis functions at the rate $N^{-\gamma}$, namely, 
\begin{align*}
    \forall N, \qquad \inf_{a_1,\ldots, a_k} \left\|\sum_{k=1}^N a_k \phi_k(x) - p_0(x)\right\|_\infty \lesssim N^{-\gamma}
\end{align*} for a sequence of $C^\infty$ functions $\{\phi_l(x)\}_{l=1}^\infty$, then there exists a neural network with number of parameters $S$ that can  approximate $s_t(x)$ at the error rate $\frac{1}{\sigma_t} S^{-\gamma}$ uniformly for all the $t$. Based on this, \cite{oko2023diffusion} establishes an approximation result that can further yield minimax optimal estimation error for $(\beta,C)$-smooth function class by adopting B-spline basis functions giving $\gamma=\beta/d$. This idea can be further extended to the Besov function class. However, this will suffer from the curse of dimensionality, resulting in a slow rate when ambient dimension $d$ is relatively large compared with the smooth index $\beta$, and the proof strategy cannot yield an optimal error bound for the class $\mathcal{H}(d, d^*, \beta, C)$ we study here.

There are also some preliminary attempts showing that neural networks can be adaptive to low-dimensional structures like $p_0 \in \mathcal{H}(d, d^*, \beta, C)$, for example, \cite{cole2024score} and \cite{kwon2025nonparametric}. However, the analysis of \cite{cole2024score} depends on the additional fundamental assumption that all the derivatives of $\log p_0(x)$ are upper bounded by a constant $c<0.5$, and the obtained rate is not optimal. \cite{kwon2025nonparametric} achieve a similar result to our Theorem \ref{theorem:approximation} and further Theorem \ref{theorem:main} when $p_0\in \mathcal{H}(d, d^*, \beta, C)$, but they adopt their own designed weight-sharing network that essentially performs convolution efficiently, and their results are inapplicable to fully-connected or even standard convolutional neural networks. 

For the class $\mathcal{H}(d, d^*, \beta, C)$, it is easy to see that the diffused function $p_t(x)$ admits approximately low-dimensional interaction structure when $t$ is small enough, i.e., $\sigma_t \lesssim (WL)^{-\beta/d^*}$, and is super-smooth akin to a $C^\infty$ function when $t$ is large, i.e., $\sigma_t \gtrsim (WL)^{-1/d}$. The main difficulty here is to establish the same fast rate in \eqref{equation:approx-error} when $(WL)^{-\beta/d^*}\ll \sigma_t \ll (WL)^{-1/d} $. This is the regime where the diffused noise $\sigma_t$ is not negligible such that the low-dimensional factorization structure no longer exists, while $\sigma_t$ is also not large enough such that $p_t(x)$ is not super-smooth. Our main argument is that $p_t(x)$ can be expressed as a composition of either low-dimensional or super smooth functions under the general regime $t\in [\underline{T}, \overline{T}]$. We offer a sketch of the proof for Theorem \ref{theorem:approximation} in Section \ref{sec:sketch}. 

\begin{remark}
One improvement compared with previous work like \cite{oko2023diffusion} and \cite{kwon2025nonparametric} is that we no longer need a technical assumption stating that $p_0$ is super smooth near the boundary of the support, i.e., Assumption 2.6 in \cite{oko2023diffusion} and Assumption B in \cite{kwon2025nonparametric}. In their proofs, such a technical assumption is imposed to establish approximation error for the score function $\nabla \log p_t(x) = \frac{\nabla p_t(x)}{p_t(x)}$ outside the region $x\in [-m_t, m_t]^d$, where $p_t$ is small thus the approximation error of $1/p_t(x)$ cannot be upper bounded at the optimal rate without this technical assumption. In our proof, we introduce a super-smooth re-weighting function $\Phi(x, t)$ such that $p_t(x)/\Phi(x, t)$ is lower bounded from below. Given super-smooth nature of $\Phi(x, t)$, we can approximate $\Phi/p_t$ and $(\Phi \sigma_t)^{-1}\nabla p_t(x)$ both at the rate of $\tilde{O}((WL)^{-2\beta/d^*})$, this allows us to get rid of such a technical assumption. 
\end{remark}

\begin{remark}
Compared with previous work \citep{oko2023diffusion, kwon2025nonparametric}, we achieve optimal approximation error for the diffused score function using deep neural networks with arbitrary (large enough) depth and width, and more importantly, do not require the use of neural networks with sparse weights. This is more practical, and requires much more effort in theory \citep{lu2021deep, kohler2021rate}. On the contrary, \cite{oko2023diffusion, kwon2025nonparametric} use a parallelization of $W$ sub-networks with width $\widetilde{O}(1)$ and depth $\widetilde{O}(1)$ for approximating score functions when $t \gtrsim (WL)^{-4/d}$, but this inevitably place sparse constraint on neural networks classes (width $\widetilde{O}(W)$ with $\widetilde{O}(W)$ non-zero parameters). 
\end{remark}

\begin{remark}
Though our approximation result is stated specifically for
densities of the functional form
$p_0 = \exp(\sum_{|J|\le d^*} f_J(x_J))$, our main arguments and the results also holds when $p_0=g(\sum_{|J|\le d^*} f_J(x_J))$ for any $C^\infty$ function $g$. 
\end{remark}

\subsection{Convergence rate in TV distance}

With the help of the above approximation error result, we are ready to state the convergence rate in TV distance between the distributions of the generated samples $\widehat{X}$ and the ground-truth $X$ when $\log p_0(x)$ admits a low-dimensional interaction structure. The following Condition \ref{condition:parameter} specifies the hyperparameters to attain the optimal rate.

\begin{condition}[Choice of parameters]\label{condition:parameter} 
Suppose $W, L$ satisfies $WL = n^{\frac{d^*}{2(2\beta+d^*)}}$ and $\{W, L\} \ge \brak{1+\log(WL)}^{c_5}$, and we pick other hyper-parameters as $\underline{T}=(WL)^{-c_3}$, $\overline{T}=c_4 \log (WL)$, $B=c_6$, where $c_3,\ldots, c_6$ are given in Theorem \ref{theorem:approximation}. 
 \end{condition}

\begin{theorem}[Error bound in TV distance]\label{theorem:main}
Assume Conditions  \ref{condition:lowerbound}, \ref{condition:betat}, and \ref{condition:parameter} hold.
There exists constants  $c_8,c_9  $ depending  only on $(c_1,c_2,\beta,d,d^*,C)$, such that for any $p_0 \in \mh(d,d^*,\beta,C)$, and $ n \ge c_8$ we have
\begin{align*}
\dse_{\{X_{0,i}\}_{i=1}^n} \left[ \mathsf{TV}(p_0,p_{\widehat{X}})\right] \leq
c_8{n^{-\frac{ \beta}{2\beta+d^*}}} (\log n)^{c_9} .
\end{align*}
Here $p_{\widehat{X}}$ is the distribution of the samples generated in Algorithm \ref{algorithm: diffusion}.
\end{theorem}

Note that $\mathcal{H}(d,d^*, \beta, C)$ contains at least all the $d^*$-variate $(\beta,C)$-smooth functions. Thus, the lower bound of the minimax optimal estimation error over $d^*$-variate $(\beta,C)$-smooth functions (\cite{khas1979lower,stone1982optimal}) gives
\begin{align*}
    \inf_{\hat{p}} \sup_{p_0 \in \mathcal{H}(d,d^*, \beta, C)} \dse_{\{X_{0,i}\}_{i=1}^n \sim p_0} \left[ \mathsf{TV}(p_0,\hat{p})\right] \gtrsim n^{-\frac{ \beta}{2\beta+d^*}},
\end{align*} which implies that the error rate in Theorem \ref{theorem:main} is minimax optimal (up to logarithmic factors). It is worth noticing that both the diffusion algorithm and the neural network architecture are blind to the interaction structure in Definition \ref{definition:interaction}. Therefore, our result indicates that the diffusion model can be adaptive to the low-dimensional interaction structure and is free of the curse-of-dimensionality $n^{-\beta/(2\beta+d)}$ rate in a structure-agnostic manner. 

\begin{remark}
It is also worth mentioning that our stochastic analysis applies to neural networks with unbounded weights \citep{farrell2021deep, kohler2021rate, fan2024noise}. This is the regime where the previous proof strategy \citep{oko2023diffusion, kwon2025nonparametric, tang2024adaptivity} cannot apply because they rely on bounded weights constraints to establish covering number in $L_\infty$ norm; see technical novelties and discussions in \aosversion{Appendix D}{Appendix \ref{section:generalerror}}.
\end{remark}

\section{Proof sketch of the neural network approximation result}
\label{sec:sketch}

To illustrate the key idea of our proof, here we provide a proof sketch for approximating the diffused density $p_t$ when $\beta=1$ and $d^*=2$. The approximation for the corresponding score function $s_t$ with $\beta\le 1$ and $d^*>2$ is similar. Additional techniques are required when $\beta>1$, which will be sketched at the end of this section.

Here to simplify the presentation, we consider using neural networks with $\widetilde{O}(1)$ depth and width $\widetilde{O}(N)$ to achieve the approximation error $\widetilde{O}(N^{-2\beta/d^*})=\widetilde{O}(N^{-1})$. We also assume $d=O(1)$. Letting $N$ be the width of the neural network considered, we consider the approximation of $p_t$ in the region $\Omega=\{(t, x): t\in [N^{-10}, N^{-1/(10d)}], x\in [-m_t, m_t]^d\}$. The approximation of $p_t$ outside the region is potentially easier. Under this time horizon of $t$, one has $m_t \in [1/2, 2]$ when $N$ is large, both $m_t$ and $\sigma_t$ are $C^\infty$ functions by Condition \ref{condition:betat}.

Denote $\mathcal{S}:=\{(i,j):i,j\in [d]\}$. When $d^*=2$ and $\beta=1$, one can write $p_0(x)=\prod_{(i,j)\in \mathcal{S}}f_{i,j}(x_i,x_j)$, where $\{f_{i,j}(x_i, x_j)\}_{i,j\in [d]}$ are all Lipschitz functions. Our analysis is based on a novel decomposition of $p_t$ that factorizes $p_t$ as a sum of terms, where each term is a product of either low-dimension or super-smooth functions. For any $i,j\in [d]$, define the the residual $\Delta_{i,j}:\real^5\rightarrow\real$ as 
\begin{align*}
    \Delta_{i,j}(t,x_i,x_j,y_i,y_j):=f_{i,j}\brak{\frac{x_i+\sigma_t y_i}{m_t},\frac{x_j+\sigma_t y_j}{m_t}}-f_{i,j}
    \brak{\frac{x_i}{m_t},\frac{x_j}{m_t}}.
\end{align*}
Let $K(y) = (2\pi)^{-d/2} \exp\brak{-\|y\|^2/2}$ be the Gaussian density.  Then, it follows from the definition of the $p_t$ and the change of variable formula that
\begin{align}
    p_t(x) &= (m_t)^{-d} \int \overbrace{\prod_{(i,j)\in \mathcal{S}} f_{i,j}\left(\frac{x_i+\sigma_t y_i}{m_t}, \frac{x_j+\sigma_t y_j}{m_t}\right)}^{p_0((x+\sigma_t y)/m_t)} K(y) dy \nonumber\\
    &\overset{(a)}{=} (m_t)^{-d} \sum_{\mathcal{A}\subset\mathcal{S}}
\underbrace{\prod_{(i,j)\in \mathcal{S}\backslash \mathcal{A}}
f_{i,j}\brak{\frac{x_i}{m_t},\frac{x_j}{m_t}}}_{G_{\mathcal{S}\backslash \mathcal{A}}(x, t)}
\underbrace{\int \prod_{(i,j)\in   \mathcal{A}}
\Delta_{i,j}K(y)dy}_{\Delta_{\mathcal{A}}(x, t)}. \label{equation:decompositionreview}
\end{align} where in $(a)$ we apply $\prod_{(i,j) \in \mathcal{S}} (a_{i,j}+b_{i,j}) = \sum_{\mathcal{A} \subset \mathcal{S}} \prod_{(i,j)\in \mathcal{A}} a_{i,j} \prod_{(i,j)\notin \mathcal{A}} b_{i,j}$ with $a_{i,j}=\Delta_{i,j}(t, x_i, x_j, y_i, y_j)$ and $b_{i,j}=f_{i,j}(x_i/m_t, y_i/m_t)$. 

The approximation of $G_{\mathcal{S}\setminus \mathcal{A}}(x,t)$ is relatively simple: it is a finite composition of several bivariate Lipschitz functions $\{f_{i,j}\}_{(i,j)\in \mathcal{S}\setminus \mathcal{A}}$ and univariate $C^{\infty}$ function $m_t$, the product function $(x,y)\rightarrow xy$ and the reciprocal function $x\rightarrow 1/x$ in the region $[1/2,2]$ given $m_t\in [1/2,2]$. Therefore, it follows from the approximation result for bivariate Lipschitz functions (see Theorem 1.1 in \cite{lu2021deep}) and the composition nature of neural networks that there exists a neural network $\widetilde{G}_{\mathcal{S}\setminus \mathcal{A}}(x,t)$ with depth $O(1)$ and width $O(N\log N)$ such that
\begin{align}
\label{equation:sketch-g}
    \forall (t, x)\in \Omega, \qquad |\widetilde{G}_{\mathcal{S}\setminus \mathcal{A}}(t, x) - {G}_{\mathcal{S}\setminus \mathcal{A}}(t, x)| \lesssim N^{-2\beta/d^*} = N^{-1}.
\end{align}

The approximation for $\Delta_{\mathcal{A}}(x, t)$ is much more involved. For given $\mathcal{A}$, let $t_*$ be the time step satisfying $\sigma_{t_*} = N^{-\frac{1}{|\mathcal{A}|}} (\log N)^7$. Given the monotonicity of $\sigma_t$, one has $\sigma_t \le \sigma_{t_*}$ if $t\le t_*$ and $\sigma_t \ge \sigma_{t_*}$ if $t\ge t_*$. Our main argument is that $\Delta_{\mathcal{A}}(x, t)$ is small enough when $t \le t_*$, and it can be decomposed into a composition of either super-smooth or low-dimensional functions when $t \ge t_*$. To see this, one has $|\Delta_{i,j}(t, x_i, x_j, y_i, y_j)| \lesssim \sigma_t$ given $f_{i,j}$ is Lipschitz function and $1/m_t \le 2$, this means
\begin{align}
\label{equation:sketch-small}
    \forall t\in [N^{-10}, t_*], \qquad |\Delta_{\mathcal{A}}(x, t)| \lesssim (\sigma_t)^{|\mathcal{A}|} 
    \leq N^{-1} (\log N)^{10d^2}.
\end{align}
On the other hand, similar to the decomposition in \eqref{equation:decompositionreview}, one also has
\begin{align}
\label{eq:decomposition-delta-large}
    \Delta_{\mathcal{A}}(x, t) = \sum_{\mathcal{B}\subset\mathcal{A}} (-1)^{\abs{\mathcal{A}\backslash\mathcal{B}}}
&\underbrace{\prod_{(i,j)\in \mathcal{A}\backslash \mathcal{B}}
 f_{i,j}\brak{\frac{x_i}{m_t},\frac{x_j}{m_t}}}_{G_{\mathcal{A}\setminus \mathcal{B}}(x, t)}\nonumber\\
&\times\underbrace{\int \prod_{(i,j)\in   \mathcal{B}}
f_{i,j}\brak{\frac{x_i+\sigma_t y_i}{m_t},\frac{x_j+\sigma_t y_j}{m_t}}K(y)dy}_{p_{t, \mathcal{B}}(x)}.
\end{align}

Note that for any fixed $t,$ $p_{t, \mathcal{B}}$ is a function dependent on at most $2|\mathcal{B}|  \le 2|\mathcal{A}| $ variables. Moreover, our 
\aosversion{Lemma E.2}
{Lemma \ref{proposition:derivativesbound}} claims that all the $\ell$-th order derivatives of $p_{t, \mathcal{B}}$ are upper bounded by $c_1(c_{\texttt{dr},1}\ell^4 )^{\ell}\sigma_t^{-\ell}$ for a constant $c_{\texttt{dr},1}$ only dependent on $(c_1,c_2,d,d^*,\beta,C)$. Now we apply 
\aosversion{Theorem A.18}{Lemma \ref{lemma:approxsmoothness}}, a finer argument of Theorem 1.1 in \cite{lu2021deep} that is uniformly in time $t$.
Then, there exists a neural network $\tilde{p}_{t, \mathcal{B}}(x, t)$ with depth $O(  \log^{16}  N )$ and width $O(N )$ such that, $\forall (t, x) \in \Omega ~\text{with}~ t\ge t_*$ we have for 
$\ell= \log(N) $,  
\begin{align*}
    |\tilde{p}_{t, \mathcal{B}}(x, t) - p_{t, \mathcal{B}}(x, t)| &\le
    {\ell}^d \norm{\nabla^{\ell}p_t}_{\infty} N^{-\frac{\ell}{ \abs{\mathcal{B}}}}
\leq c_1  \ell^d \brak{\frac{c_{\texttt{dr},1} \ell^4}{\sigma_{t} N^{\frac{1}{ \abs{\mathcal{B}}}}}}^\ell\\
    &\overset{(a)}{\le} c_1(\log N)^d \brak{\frac{\log N}{c_{\texttt{dr},1}}}^{-\log N} = O(N^{-10}),
\end{align*} where $(a)$ follows from the fact that $\abs{\mb}\leq \abs{\ma}$ and the time step is large $t \ge t_*$ such that $\sigma_t N^{1/|\mathcal{A}|} \ge \log^5(N) \ge \log(N) \ell^4$. 
Therefore, for any $\mathcal{A}$, we can use a neural network $\widetilde{\Delta}_{\mathcal{A}}(x, t)$ with depth $O(\log^2 N)$ width $O(N\log^{2d} N)$ to approximate $1\{t \ge t_*\} \times \sum_{\mathcal{B}\subset\mathcal{A}} \widetilde{G}_{\mathcal{A}\setminus \mathcal{B}}(x, t) \cdot \tilde{p}_{t, \mathcal{B}}(x, t)$ with error $O(N^{-10})$. Combining it with \eqref{equation:sketch-small}, $\widetilde{\Delta}_{\mathcal{A}}(x, t)$ satisfies
\begin{align}
\label{equation:sketch-delta}
    \forall (t, x)\in \Omega, \qquad |\widetilde{\Delta}_{\mathcal{A}}(x, t) - \Delta_{\mathcal{A}}(x, t)| \lesssim N^{-1} (\log N)^{10d^2}.
\end{align}

Now \eqref{equation:sketch-g} and \eqref{equation:sketch-delta} together assert that both $G_{\mathcal{S}\setminus \mathcal{A}}(x, t)$ and $\Delta_{\mathcal{A}}(x, t)$ can be approximated by neural network with depth $\widetilde{O}(1)$ and width $\widetilde{O}(N)$ at the error rate $\widetilde{O}(N^{-1})$. Combing it with the decomposition \eqref{equation:decompositionreview} at the beginning, which states the target $p_t(x, t)$ is a sum of products of these terms $\{G_{\mathcal{S}\setminus \mathcal{A}}(x, t), \Delta_{\mathcal{A}}(x, t)\}_{\mathcal{A}\subseteq \mathcal{S}}$, we can find a neural network $\tilde{p}(t, x)$ with depth $\widetilde{O}(1)$ and width $\widetilde{O}(N)$ that can approximate $p_t(x, t)$ uniformly well at the error rate $\widetilde{O}(N^{-1})$.

\subsection{Extension to $\beta>1$}

It is worth noticing that the current decomposition \eqref{equation:decompositionreview} cannot yield optimal approximation error when $\beta>1$. In this part, we illustrate the potential issue and a finer approximation when $\beta=2$; the extension to general $\beta>1$ is similar.

Turning to the decomposition \eqref{equation:decompositionreview}, similar to the discussion above, for any $\mathcal{A} \subseteq \mathcal{S}$, there exists a neural network $\tilde{G}_{\mathcal{S}\setminus\mathcal{A}}(x, t)$ with depth ${O}(1)$ and width $O(N\log N)$ such that
\begin{align*}
    \forall (t, x)\in \Omega, \qquad |\widetilde{G}_{\mathcal{S}\setminus \mathcal{A}}(t, x) - {G}_{\mathcal{S}\setminus \mathcal{A}}(t, x)| \lesssim N^{-2\beta/d^*} = N^{-2}.
\end{align*}


However, the arguments for $\Delta_{\mathcal{A}}(x, t)$ above cannot give the optimal approximation error for it. To be specific, following a similar discussion to the decomposition in \eqref{eq:decomposition-delta-large}, we can find a neural network $\widetilde{\Delta}_{\mathcal{A}}(x, t)$ with depth $\widetilde{O}(1)$ and width $\widetilde{O}(N)$ such that
\begin{align*}
    \forall (t, x)\in \Omega \cap \{(t, x): t\geq t_{*,1} \}, \qquad |\widetilde{\Delta}_{\mathcal{A}}(x, t) - \Delta_{\mathcal{A}}(x, t)| = \widetilde{O}(N^{-2}).
\end{align*} with $t_{*,1}$ is the timestep satisfying $\sigma_{t_{*, 1}} = N^{-\frac{1}{|\mathcal{A}|}} (\log N)^7$. Here the choice of $t_{*,1}$ is to let $\sigma_t N^{1/|\mathcal{B}|} \ge \log^7(N)$ for any $t\ge t_{*, 1}$ in the approximation of $p_{t, \mathcal{B}}(x)$ when $\mb\subset \ma$. The main issue here is that we no longer have $|\Delta_{\mathcal{A}}(x, t)| = \widetilde{O}(N^{-2})$ when $t \le t_{*, 1}$. This is because akin to \eqref{equation:sketch-small}, one only has
\begin{align*}
    |\Delta_{\mathcal{A}}(x, t)| \le (\sigma_t)^{|\mathcal{A}| (\beta \land 1)} = \sigma_t^{|\mathcal{A}|}.
\end{align*} This implies that $|\Delta_{\mathcal{A}}(x, t)| = \widetilde{O}(N^{-2})$ only applies when $\sigma_t = \widetilde{O}( N^{-\frac{2}{|\mathcal{A}|}} ) = o(\sigma_{t_{*, 1}})$. 

 \begin{figure}
     \centering
     \subfigure[Approximator for $\beta\leq 1$]{\includegraphics[width=0.4\linewidth]{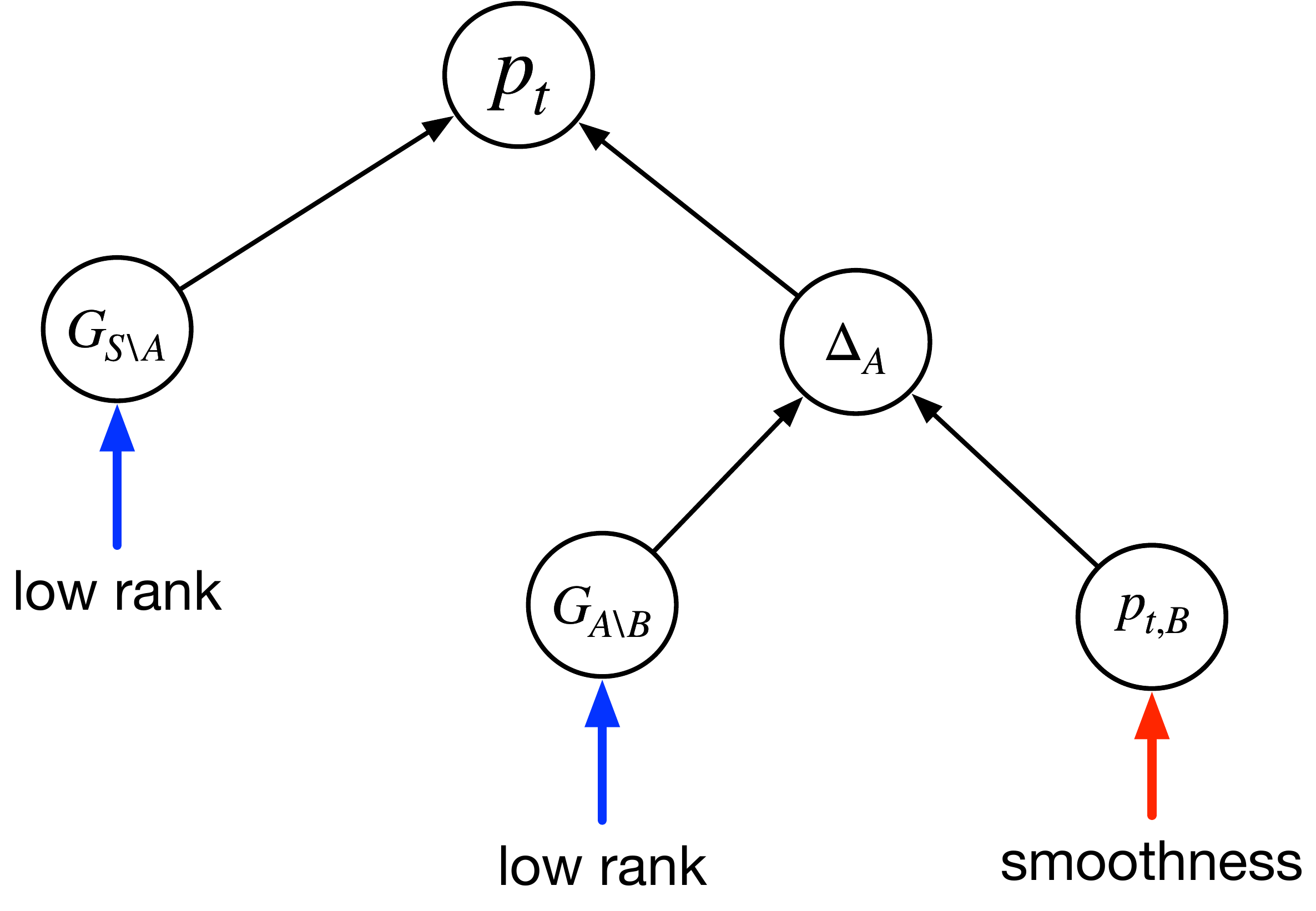}}
     \subfigure[Approximator for $\beta> 1$]{\includegraphics[width=0.4\linewidth]{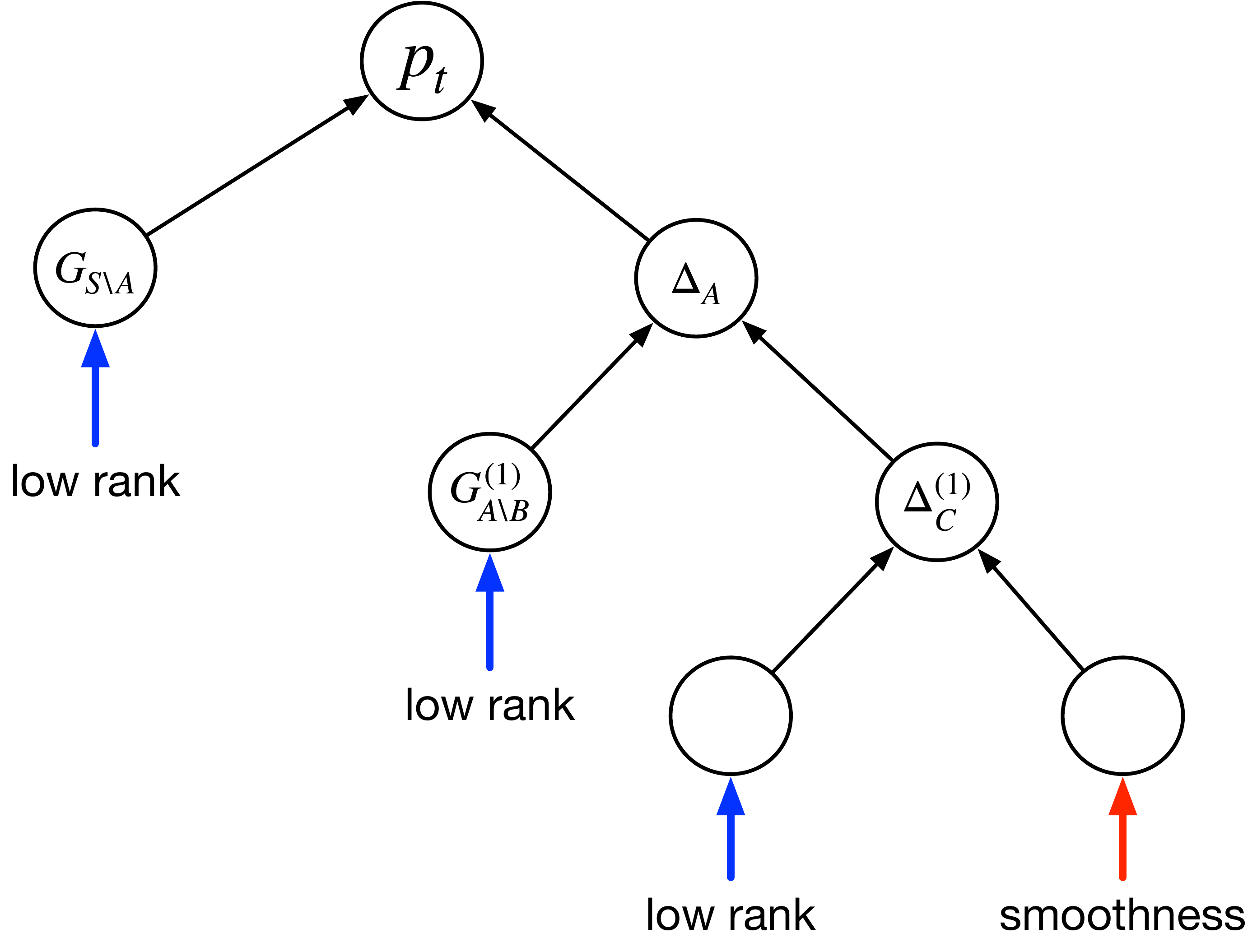}}
    \caption{Figures illustrating the construction of the neural network approximator for $p_t$.  
    The approximator for $p_t$ is constructed in a bottom-to-top manner.
    The blue arrow  from ``low rank'' to node $x$ indicates that $x$ is a product of  $d^*$-dimensional H{\"o}lder-smooth functions, whose neural network approximators are constructed with 
    Theorem 1.1 in \cite{lu2021deep}. For example, the nodes $G_{\mathcal{S}\backslash\mathcal{A}}, G_{\mathcal{A}\backslash\mathcal{B}}$   are product of  $d^*$-dimensional $(\beta,C)$-functions, and 
    the node $G^{(1)}_{\mathcal{A}\backslash\mathcal{B}}$ is a product of  $d^*$-dimensional $(\beta-1,C)$-functions.
    The red arrow from ``smoothness'' to node $x$ indicates that $x$ is super-smooth or its parent node is negligible, whose neural network approximator is constructed with \aosversion{Theorem A.18}{Lemma \ref{lemma:approxsmoothness}}.
    For example, $p_{t,\mathcal{B}}$ is super-smooth otherwise $\Delta_{\mathcal{A}}$ is negligible.
    Black arrows go from nodes $x,y$ to $z$ indicate that $z$ is a sum-of-product of $x$ and $y.$ Thus the neural network approximator for $z$ can be constructed using the approximators for $x$ and $y.$  
    }
     \label{fig:proofsketch}
 \end{figure}

Therefore, a finer approximation of $\Delta_{\mathcal{A}}(x, t)$ in the time interval $t\in [0, t_{*, 1}]$, especially when $\sigma_t \gg N^{-\frac{2}{|\mathcal{A}|}}$, is required, this is also where high-order smoothness should come into play in the decomposition. 
We first introduce some additional notations. Let $\tilde{x}_j = x_j / m_t$ and $\tilde{\sigma}_t = \sigma_t / m_t$. We let
$D_{i,j}(\tilde{x}_i, \tilde{x}_j) = [\frac{\partial f_{i,j}}{\partial x_i}(\tilde{x}_i, \tilde{x}_j), \frac{\partial f_{i,j}}{\partial x_j}(\tilde{x}_i, \tilde{x}_j)]^\top \in \mathds{R}^2$ and define the scaled second order Taylor residual as
\begin{align*}
    \Delta_{i,j}^{(1)}(t, x_i, x_j, y_i, y_j) = \frac{1}{\tilde{\sigma}_t} \left\{\Delta_{i,j}(t, x_i, x_j, y_i, y_j) - \tilde{\sigma}_t \cdot [y_i, y_j] D_{i,j} (\tilde{x}_i, \tilde{x}_j) \right\}.
\end{align*}
Plugging the identity $\Delta_{i,j}(t, x_i, x_j, y_i, y_j)=\tilde{\sigma}_t \{\Delta_{i,j}^{(1)}(t, x_i, x_j, y_i, y_j) + [y_i, y_j] D_{i,j} (\tilde{x}_i, \tilde{x}_j)\}$ into the definition of $\Delta_{\mathcal{A}}(x, t)$ in \eqref{equation:decompositionreview} then gives that
\begin{align*}
    \Delta_{\mathcal{A}}(x, t) 
    &= (\tilde{\sigma}_t)^{|\mathcal{A}|} \sum_{\mathcal{C} \subseteq \mathcal{B} \subseteq \mathcal{A}} \underbrace{\prod_{(i,j)\in \mathcal{B}\setminus \mathcal{C}} D_{i,j,2}(\tilde{x}_i,\tilde{x}_j)}_{G^{(1)}_{\mathcal{B}\setminus \mathcal{C},2}(x, t)} \times \underbrace{\prod_{(i,j)\in \mathcal{A}\setminus \mathcal{B}} D_{i,j,1}(\tilde{x}_i,\tilde{x}_j)}_{G^{(1)}_{\mathcal{A}\setminus \mathcal{B},1}(x, t)} \\
    &~~~~~~~~~~~~~~~~~~~~~~~~~~~~ \times \underbrace{\int \prod_{(i,j)\in \mathcal{C}} \Delta_{i,j}^{(1)} \prod_{(i,j)\in \mathcal{A}\setminus \mathcal{B}} y_i \prod_{(i,j)\in \mathcal{B}\setminus \mathcal{C}} y_j K(y) dy}_{\Delta^{(1)}_{\mathcal{C}}(x, t)}.
\end{align*}
The idea is, the functions in $\{G_{\mathcal{K},l}^{(1)}(x, t)\}_{\mathcal{K}\subseteq \mathcal{S}, l\in [2]}$ can be approximated by neural networks at the error rate $\widetilde{O}(N^{-(\beta-1)})=\widetilde{O}(N^{-1})$ given each function itself is a composition of simple functions and $2$-variate, $(1, C)$ smooth functions. On the other hand, the same approximation error rate $\widetilde{O}(N^{-(\beta-1)})=\widetilde{O}(N^{-1})$ is also attainable for $\Delta_{\mathcal{C}}^{(1)}(x, t)$.
The idea is similar to $\Delta_{\mathcal{A}}(x, t)$ in the previous discussion: $\Delta_{\mathcal{C}}^{(1)}(x, t)$ is negligible, i.e., $|\Delta_{\mathcal{C}}^{(1)}(x, t)| \le \widetilde{O}(N^{-1})$, when $\sigma_t \le N^{-\frac{1}{|\mathcal{C}|}}  (\log N)^7$, and can be decomposition into a composition of either at least $2$-variate $(1, C)$-smooth functions, or super-smooth functions similar to \eqref{eq:decomposition-delta-large} as long as $\sigma_t \ge N^{-\frac{1}{|\mathcal{C}|}} (\log N)^7$.

Therefore, when $t \le t_{*, 1}$, one can construct a neural network that approximates $\Delta_{\mathcal{A}}(x, t)$ at the rate of
\begin{align*}
    \tilde{\sigma}_t^{|\mathcal{A}|} \times \widetilde{O}(N^{-1}) \lesssim (\sigma_{t_{*, 1}})^{|\mathcal{A}|} \times \widetilde{O}(N^{-1}) = \widetilde{O}(N^{-2}).
\end{align*}
The proof idea is summarized in the right panel of Figure \ref{fig:proofsketch}.

\section{Convergence rates in $W_1$ distance}\label{section: w1}
Besides the TV distance, we also measure the performance of the diffusion model distribution estimator in the Wasserstein distance of order one ($W_1$ distance). 
For any distributions $X$ and $\widehat{X}$, 
let $\Pi(X, \widehat{X})$ denote all joint distributions $\pi$  that have marginals $X$ and $\widehat{X}$.
The  $W_1$ distance between $X$ and $\widehat{X}$ is defined as
$$
W_1(X, \widehat{X}):= \inf _{\pi \in \Pi(X, \widehat{X})} \int\|x-y\|  d \pi(x, y) = \sup_{f: \mbox{ \scriptsize Lip}(f) \leq 1} |E f(X) - f(\widehat{X})|.
$$

\medskip
\noindent
\textbf{Piecewise score estimator.}
To achieve an efficient estimator in $W_1$ distance, we slightly modify the 
score-matching procedure. We partition the time horizon $[\ut,\ot]$ as consecutive intervals 
$t_0=\ut<t_1<\cdots <t_P=\ot$  for time parameters specified later. The piecewise score estimator uses neural networks with varying sizes in time intervals $[t_{\ell-1}, t_{\ell}]$ for $\ell\in [P]$.
Denote the score-matching loss on the time interval $[t_{j-1},t_j]$ as
\begin{align}\label{equation:lossinterval}
    \widehat{\mathsf{L}}_{n,j}(s) = \frac{1}{n} \sum_{i=1}^n \dse_{t\sim \mathsf{Unif}(t_{j-1},t_j), X_t\sim p_{t|0}(\cdot|X_{0,i})} \left[ \left\|s(X_t, t)-\nabla \log p_{t|0}(X_t | X_{0, i})\right\|_2^2\right].
\end{align} 
where $\mathsf{Unif}([t_{j-1},t_j])$ is the uniform distribution  on the time interval $[t_{j-1},t_j]$. 
We then train the neural network score-matching estimator on $[t_{j-1},t_j]$ as
\begin{align}\label{equation:estimatorinterval}
\hat{s}_j:=\arg\min_{s\in \mathcal{F}_{\mathsf{NN}}(d+1,L_j,W_j,B_j)}
\widehat{\mathsf{L}}_{n,j}(s),
\end{align} 
for the network parameters $(L_j,W_j,B_j)$ specified later. 
The final score estimator $\hat{s}$ is defined as 
\begin{align}\label{equation:estimatorw1}
   \hat{s}(t,x)=\hat{s}_j(t,x)\text{ for any } t\in [t_{j-1},t_j],\, j\in [P].
\end{align}
The underlying reason behind adopting the piecewise score estimator is that the score function estimation error in different timesteps contributed to the final $W_1$ error in a non-uniform manner. In particular, the score function estimation error for large timesteps contributes less to the final $W_1$ error, hence using ReLU neural network score function estimators with descending width/depth in timestep as in Condition \ref{condition:parameterw1} can attain a slightly faster rate; see also \cite{oko2023diffusion}.

\medskip
 \noindent \textbf{Diffusion model with piecewise score-matching.} 
We propose the distribution estimator $\widehat{X}^{\texttt{W}}$ for $X$ by solving the backward process 
with the estimated score function $\hat{s}$ in \eqref{equation:estimatorw1}. 
The process is summarized in Algorithm \ref{algorithm: diffusionW1}.  
\begin{algorithm}[!ht]
    \renewcommand{\algorithmicrequire}{\textbf{Input:}}
    \renewcommand{\algorithmicensure}{\textbf{Output:}}
    \caption{Diffusion models with piecewise score-matching.}
    \label{algorithm: diffusionW1}
    \begin{algorithmic}[1]
        \Require i.i.d. samples $\{X_{0,i}\}_{i=1}^n$ drawn from the distribution $p_0$.
        \Ensure an output $\widehat{X}^{\texttt{W}}$ whose distribution is close to $p_0$.
        \State Initialization: $\hat{s}\equiv0$.
        \For{$j=1$ to $P$} 
        \State 
        Set the network parameters (depth $L_j$, width $W_j$, and truncation threshold $B_j$) and the
        time parameter $t_j$
          as specified in Condition \ref{condition:parameterw1}. 
          \State 
          Compute the estimated score function $\hat{s}_j$ on $[t_{j-1},t_j]$
         by minimizing the loss  $ \widehat{\mathsf{L}}_{n,j}$   defined in \eqref{equation:lossinterval} as
        $
         \hat{s}_j:=\arg\min_{s\in \mathcal{F}_{\mathsf{NN}}(d+1,L_j,W_j,B_j)} \widehat{\mathsf{L}}_{n,j}(s).
         $
         \State Update the estimated score function
        $\hat{s}(t,x)\leftarrow\hat{s}(t,x)+\hat{s}_j(t,x)\mathds{1}\{t_{j-1}\leq t \leq t_j\}.$
        \EndFor
        \State Draw the sample from the backward process with the estimated score $\hat{s}$:
        \begin{align}\label{equation:backwardestimatedw1}
        & d \widehat{Y}_t=\beta_{\overline{T}-t}\left(\widehat{Y}_t+2 \hat{s}(\widehat{Y}_t, \overline{T}-t)\right) d t+
        \sqrt{2 \beta_{\overline{T}-t}} d B_t, \, \widehat{Y}_0 \sim N(0,I_d).
        \end{align}
    \State \Return $\widehat{X}^{\texttt{W}}$ as $\widehat{Y}_{\overline{T}-\underline{T}}$.
 \end{algorithmic}
\end{algorithm}

\subsection{Convergence rate in $W_1$ distance}
We state the convergence rate in $W_1$ distance between the generated samples $\widehat{X}^{\texttt{W}}$ and the ground-truth $X$ when $\log p_0(x)$ admits a low-dimensional interaction structure. The following Condition \ref{condition:parameterw1} specifies the hyperparameters to attain the desired rate.

\begin{condition}[Choice of parameters]\label{condition:parameterw1} 
Suppose $P=\lfloor\log_2(\ot/\ut)\rfloor+1$, $t_0=\ut$, $t_1=n^{-\frac{2d^*}{d(2\beta+d^*)}}$, $t_{j+1}=2t_j \wedge \ot$ for any $j\in [P-1]$.
Suppose $W_j$ and $L_j$ satisfy $W_jL_j =t_j^{-\frac{d}{4}}(\log W_1L_1)^{2c_5+500d} $ and 
$W_j\wedge L_j\geq (1+\log (W_1L_1))^{c_5}$
for any  $j\in [P]$, and we pick other hyper-parameters as $\underline{T}=(W_1L_1)^{-c_3}$, $\overline{T}=c_4 \log (W_1L_1)$, $B_j=c_6$, where $c_3,\ldots, c_5$ are given in Theorem \ref{theorem:approximation}.
 \end{condition}

\begin{theorem}[Error bound in $W_1$ distance]\label{theorem:mainw1}
Assume Conditions  \ref{condition:lowerbound}, \ref{condition:betat}, and \ref{condition:parameterw1} hold.
There exists constants  $c_8,c_9  $ depending  only on $(c_1,c_2,\beta,d,d^*,C)$, such that for any $p_0 \in \mh(d,d^*,\beta,C)$, and $ n \ge c_8$ we have
\begin{align*}
\dse_{\{X_{0,i}\}_{i=1}^n} \left[ W_1(p_0,p_{\widehat{X}^{\texttt{W}}})\right] \leq 
\begin{cases}
     c_8{n^{-\frac{ \beta+\frac{d^*}{d}}{2\beta+d^*}}} (\log n)^{c_9}, &\text{ for }d\geq 2, \\\\
     c_8 n^{-1/2}(\log n)^{c_9},& \text{ for }d=1.
\end{cases}
\end{align*}
Here $\widehat{X}^{\texttt{W}}$ is the distribution of the samples generated in Algorithm \ref{algorithm: diffusionW1}.
\end{theorem}

It is known from \cite{oko2023diffusion} that diffusion model can estimate $(\beta,C)$-smooth densities with error rate $n^{-(\beta+1)/(2\beta+d)}$ in $W_1$ distance. Theorem \ref{theorem:mainw1} establishes a rate faster than that and is free of the curse of dimensionality measured in $W_1$ distance. Notably, the diffusion model achieves a strictly faster rate in $W_1$ distance than TV distance (Theorem \ref{theorem:main}).

\begin{remark}
We point out here that the sharpness of Theorem \ref{theorem:mainw1} remains unclear. Combining the minimax lower bound for estimating $d^*$-covariate $(\beta,C)$-smooth densities in $W_1$ distance \citep{liang2021well} and the simple fact that $\mathcal{H}(d^*,d^*, \beta, C) \subseteq \mathcal{H}(d,d^*, \beta, C)$ gives the lower bound
\begin{align*}
    \inf_{\hat{p}} \sup_{p_0 \in \mathcal{H}(d,d^*, \beta, C)} \dse_{\{X_{0,i}\}_{i=1}^n \sim p_0} \left[ W_1(p_0,\hat{p})\right] \gtrsim n^{-\frac{ \beta+1}{2\beta+d^*}} + n^{-1/2}.
\end{align*} There is a gap between the upper bound presented in Theorem \ref{theorem:mainw1} and the above lower bound when $d>d^*$. Exploring the matching minimax optimal estimation error in this setup and whether the diffusion model can attain such an optimality calls for new techniques, and we leave it for future studies.
\end{remark}

\bibliographystyle{apalike}
\bibliography{arxiv_main}

\end{document}